\documentclass[11pt]{article}

\usepackage{amscd,amsmath, amssymb, fancyhdr, graphicx, url}

\usepackage[backref=page]{hyperref}

\renewcommand*{\backrefalt}[4]{%
	\ifcase #1 (Not cited.)%
	\or        (Cited on page~#2.)%
	\else      (Cited on pages~#2.)%
	\fi}

\hypersetup{
	colorlinks   = true,
	citecolor    = magenta,
	linkcolor    =blue,
	urlcolor     =magenta	
}

\numberwithin{equation}{section}

\newcommand{\version}{version 4.0,\ \ June 12, 2023, journal version}

\makeatletter
\def\x@arrow{\DOTSB\Relbar}
\def\xlongrightarrowfill@{\arrowfill@\relbar\relbar\longrightarrow}
\newcommand{\xlongrightarrow}[2][]{%
        \ext@arrow 0099\xlongrightarrowfill@{#1}{#2}}
\makeatother

\def\eqref#1{(\ref{#1})}
\newcommand{\goth}{\mathfrak}

\newcommand{\arrow}{{\:\longrightarrow\:}}
\newcommand{\Z}{{\Bbb Z}}
\def\C{{\Bbb C}}
\def\P{{\Bbb P}}

\newcommand{\R}{{\Bbb R}}
\newcommand{\Q}{{\Bbb Q}}
\renewcommand{\H}{{\Bbb H}}
\newcommand{\6}{\partial}
\def\1{\sqrt{-1}\:}
\newcommand{\restrict}[1]{{\left|_{{\phantom{|}\!\!}_{#1}}\right.}}
\newcommand{\cntrct}                
{\hspace{2pt}\raisebox{1pt}{\text{$\lrcorner$}}\hspace{2pt}}

\newcommand{\cac}{{\cal C}}

\renewcommand{\tilde}{\widetilde}
\renewcommand{\bar}{\overline}
\renewcommand{\phi}{\varphi}
\renewcommand{\epsilon}{\varepsilon}
\renewcommand{\geq}{\geqslant}

\renewcommand{\max}{{\rm max}}

\newcommand{\NS}{\operatorname{\sf NS}}

\newcommand{\Pos}{\operatorname{Pos}}

\newcommand{\Aut}{\operatorname{Aut}}

\newcommand{\Mon}{\operatorname{\sf Mon}}

\newcommand{\Alb}{\operatorname{Alb}}

\newcommand{\St}{\operatorname{St}}

\newcommand{\rk}{\operatorname{rk}}

\newcommand{\Tr}{\operatorname{Tr}}

\renewcommand{\Re}{\operatorname{Re}}
\renewcommand{\Im}{\operatorname{Im}}

\newcommand{\proof}{\noindent{\bf Proof:\ }}
\newcommand{\pstep}{{\bf Proof. Step 1:\ }}

\newcounter{Mycounter}[section]
\newcounter{lemma}[section]
\renewcommand{\thelemma}{{Lemma \thesection.\arabic{lemma}}}
\newcommand{\lemma}{%
    \setcounter{lemma}{\value{Mycounter}}
    \refstepcounter{lemma}
    \stepcounter{Mycounter}
    {\noindent \bf \thelemma:\ }}

\newcounter{claim}[section]
\renewcommand{\theclaim}{{Claim \thesection.\arabic{claim}}}
\newcommand{\claim}{%
    \setcounter{claim}{\value{Mycounter}}
    \refstepcounter{claim}
    \stepcounter{Mycounter}
    {\noindent \bf \theclaim:\ }}

\newcounter{sublemma}[section]
\newcounter{corollary}[section]
\renewcommand{\thecorollary}{{Corollary \thesection.\arabic{corollary}}}
\newcommand{\corollary}{%
    \setcounter{corollary}{\value{Mycounter}}
    \refstepcounter{corollary}
    \stepcounter{Mycounter}
    {\noindent \bf \thecorollary:\ }}

\newcounter{theorem}[section]
\renewcommand{\thetheorem}{{Theorem \thesection.\arabic{theorem}}}
\newcommand{\theorem}{%
    \setcounter{theorem}{\value{Mycounter}}
    \refstepcounter{theorem}
    \stepcounter{Mycounter}
    {\noindent \bf \thetheorem:\ }}

\newcounter{conjecture}[section]
\renewcommand{\theconjecture}{{Conjecture \thesection.\arabic{conjecture}}}
\newcommand{\conjecture}{%
    \setcounter{conjecture}{\value{Mycounter}}
    \refstepcounter{conjecture}
    \stepcounter{Mycounter}
    {\noindent \bf \theconjecture:\ }}

\newcounter{proposition}[section]
\renewcommand{\theproposition}
      {{Proposition \thesection.\arabic{proposition}}}
\newcommand{\proposition}{%
    \setcounter{proposition}{\value{Mycounter}}
    \refstepcounter{proposition}
    \stepcounter{Mycounter}
    {\noindent \bf \theproposition:\ }}

\newcounter{definition}[section]
\renewcommand{\thedefinition}
      {{Definition~\thesection.\arabic{definition}}}
\newcommand{\definition}{%
    \setcounter{definition}{\value{Mycounter}}
    \refstepcounter{definition}
    \stepcounter{Mycounter}
    {\noindent \bf \thedefinition:\ }}

\newcounter{example}[section]
\renewcommand{\theexample}{{Example \thesection.\arabic{example}}}
\newcommand{\example}{%
    \setcounter{example}{\value{Mycounter}}
    \refstepcounter{example}
    \stepcounter{Mycounter}
    {\noindent \bf \theexample:\ }}

\newcounter{remark}[section]
\renewcommand{\theremark}{{Remark \thesection.\arabic{remark}}}
\newcommand{\remark}{%
    \setcounter{remark}{\value{Mycounter}}
    \refstepcounter{remark}
    \stepcounter{Mycounter}
    {\noindent \bf \theremark:\ }}

\newcounter{problem}[section]
\newcounter{question}[section]
\makeatletter

\@addtoreset{equation}{section} \@addtoreset{footnote}{section}
\makeatother

\def\blacksquare{\hbox{\vrule width 5pt height 5pt depth 0pt}}
\def\endproof{\blacksquare}

\begin{document}

\begin{center}
{\LARGE\bf
Parabolic automorphisms \\[3mm] of hyperk\"ahler manifolds\\[3mm]
}

Ekaterina Amerik\footnote{Both authors acknowledge support of 
HSE University basic research program; also partially supported by ANR (France) project FANOHK} 
Misha Verbitsky\footnote{Partially supported by 
FAPERJ E-26/202.912/2018 and CNPq - Process 313608/2017-2.

{\bf Keywords:} hyperk\"ahler manifold, ergodic theory,

{\bf 2010 Mathematics Subject
Classification:} }

\end{center}

{\small \hspace{0.15\linewidth}
\begin{minipage}[t]{0.7\linewidth}
{\bf Abstract} \\
A parabolic automorphism of a hyperk\"ahler manifold
is a holomorphic automorphism acting on $H^2(M)$ by
a non-semisimple quasi-unipotent linear map. We prove that
a parabolic automorphism which preserves a Lagrangian fibration acts on almost all fibers ergodically. 
	The existence of an invariant Lagrangian
fibration is automatic for manifolds satisfying the
hyperk\"ahler SYZ conjecture; 
this includes all known examples of hyperk\"ahler manifolds.
When there are two parabolic automorphisms preserving two
distinct Lagrangian fibrations, it follows that the group
they generate acts on $M$ ergodically. Our results
generalize those obtained by S. Cantat for K3 surfaces.
\end{minipage}
}

\tableofcontents


\section{Introduction}


\subsection{K3 surfaces}

Let $S$ be a projective K3 surface and let $f: S\to S$ be
an automorphism. The Neron-Severi lattice 
 $\NS(S)$ is of signature $(1, \rho_S-1)$ where $\rho_S$
is the Picard number of $S$, and $f$ induces an
automorphism $f^*$ of $NS(S)$. When $\rho_S\geq 3$, the
projectivization of the cone of vectors of positive square in $V=\NS(S)\otimes
\R$ can be viewed as a model of the hyperbolic space which
we denote by $\H_S$. The projectivization of the isotropic cone is the boundary of the hyperbolic space, called the absolute,
and $f^*\in O(V)$ induces an isometry
of $\H_S$. Replacing $f$ by $f^2$, we may assume $f^*\in SO^+(V)$\footnote{Here $SO^+(V)$ denotes the connected component of the unity of $SO(V)$.}. The following theorem is well-known (see \cite{_Kapovich:hyperbolic_} or \cite{R} for details).

\hfill

{\bf Theorem-definition:}
Let $n>1$, and $\alpha \in SO^+(1,n)$ a non-trivial isometry
acting on $(V, q)$, $q$ quadratic form of signature $(1,n)$. Then one and only one of these three cases occurs:
\begin{description}
\item[(i)] 
$\alpha$ has an eigenvector $x$ with $q(x)>0$ 
($\alpha$ is an {\bf elliptic isometry)} 
\item[(ii)]
	$\alpha$ has two eigenvectors $x, y$ with $q(x)=q(y)=0$ 
		and  real eigenvalues $\lambda_x$, $\lambda_y=\lambda_x^{-1}$ satisfying $|\lambda_x|>1$
($\alpha$ is a {\bf loxodromic, or hyperbolic, isometry)}
\item[(iii)] $\alpha$ has a single (up to a constant) eigenvector 
$x$ with $q(x)=0$. In this case
the eigenvalue   $\lambda_x=1$,
and $\alpha$ is called a {\bf parabolic isometry}).
\end{description}

In terms of the action on the hyperbolic space $\H$, an elliptic isometry has a fixed point in $\H$, a parabolic one has one fixed point on the absolute, and a loxodromic one has two fixed points on the absolute (and no fixed points inside $\H$). 

\hfill

\remark  All eigenvalues
of elliptic and parabolic 
isometries have absolute value 1. 
Hyperbolic and elliptic isometries
are semisimple (that is, the corresponding linear operators
are  diagonalizable over $\C$), parabolic are
not. When $(V, q)$ has some underlying integral structure, for
example if $V=\NS(S)\otimes \R$, and $q$ is the intersection form, 
the elliptic isometries preserving the integral structure
are of finite order.
 Parabolic ones are of
infinite order; indeed, any linear homomorphism of
finite order is semisimple.

\hfill

We call the automorphism $f$ {\bf parabolic} if $f^*$ is parabolic.
In this case, $f^*$ preserves a class of
self-intersection zero, which one may clearly assume nef
(by invariance of the nef cone) and integral (as $f^*$ is
integral).

Remark that the classification of automorphisms into elliptic, parabolic and loxodromic types also makes sense for compact K\"ahler surfaces, when one
replaces $\NS(S)\otimes \R$ by $H^{1,1}_{\R}(X)$. Though the latter space has no natural integral structure, it is easy to see that elliptic isometries are of finite order\footnote{And so, by Fujiki-Liebermann theorem, are the elliptic automorphisms
when $S$ has no vector fields, in particular in the K3 case.} and that the parabolic isometries fix a nef integral class, using the integral structure on $H^2(X,\R)$ (see for example \cite{M}, p. 31).

By Gizatullin's theorem (\cite{_Gizatullin_,_Grivaux_}), 
any parabolic automorphism $p$
of a projective surface admits an invariant elliptic fibration,
or ``elliptic pencil'', that is, a fibration $\pi: S\to \P^1$ with
general fiber a smooth curve of genus 1 (we do not assume that this fibration has a section), and such that $p$ sends fibers to fibers. 
For a K3 surface this is easy to see directly, indeed, one deduces from the Riemann-Roch theorem that 
the linear system of the sections of a nef line bundle of square zero is an elliptic pencil.
Moreover, by topological reasons $\pi$ has at least three singular fibers, which must be permuted by the action of $p$.
So $p$ is of 
finite order on the base, hence a power of $p$ preserves the fibers\footnote{In what follows, we say that $p$ {\bf preserves a fibration} $\pi$ when $\pi \circ p=\pi$. When $p$ just sends fibers to fibers, we call $\pi$ {\bf an invariant fibration.}}. 

These observations have been made by
Serge Cantat in \cite{_Cantat:dyn_K3_TG_}. In Proposition 2.2, Proposition 2.1 of \cite{_Cantat:dyn_K3_TG_}, he has shown that
a parabolic automorphism of a projective K3 surface has dense orbits on almost all genus 1 
fibers. Moreover (\cite[Corollaire 3.2]{_Cantat:dyn_K3_TG_}) if two automorphisms $f_1, f_2$ preserve
distinct elliptic fibrations $\pi_1, \pi_2$, then the
group they generate acts on $S$ ergodically, that is,
every measurable invariant subset has zero or full
measure\footnote{Here the measure is the one coming from the volume form $\Omega\wedge\bar{\Omega}$, where $\Omega$ is the holomorphic 2-form on $S$, unique up to a constant}. The ergodicity statement follows from the density
of orbits by Fubini theorem, so only the density statement
requires work to be done.

Our aim is to generalize these results to hyperk\"ahler
manifolds of arbitrary dimension (not necessarily projective).

Let us recall some more remarks
about elliptic pencils on K3 surfaces
and parabolic automorphisms
 from \cite{_Cantat:dyn_K3_TG_},
some of them well-known.

Let $S$ be a K3 surface. If the
intersection form on $\NS(S)$ represents zero and does not represent $-2$, then $S$
carries an elliptic pencil $\pi$. If, moreover,
the rank of $\NS(S)$ is greater than 
two, 
then there is always a parabolic automorphism 
preserving $\pi$. Indeed,\footnote{For any K3
surface $S$, the group generated by the reflexions in
the $(-2)$-curves and the image of $\Aut(S)$ is a lattice in
$O(\NS\otimes_\Z \R)$ (\cite{_PSh-Sh:Torelli_}).} $\Aut(S)$
is a lattice, hence the stabilizer of an 
isotropic class has infinite order: 
otherwise the quotient of $\H_S$ by $\Aut(S)$ would have
an end of infinite volume.

As an example of a surface admitting two parabolic
automorphisms with distinct nef eigenvectors we
 may take a sufficiently general complete intersection of type $(2,2,2)$ in $(\P^1)^3$. Indeed each of the three projections to $\P^1$ is an elliptic pencil 
and each of the three projections to $(\P^1)^2$ is a double covering. The product of any two of the three covering involutions is 
parabolic (see e.g. \cite{_Cantat:dyn_K3_TG_}, Exemple 5, and \cite{W} for details).

\subsection{Holomorphic symplectic manifolds}

The purpose of the present article is to generalize this
picture to higher dimension, that is, for arbitrary
 {\bf irreducible holomorphic symplectic
  manifolds}, also called {\bf maximal holonomy
  hyperk\"ahler manifolds}.\footnote{We use both names interchangeably, so that irreducible holomorphic symplectic manifolds are always assumed compact and K\"ahler.} The role of the
intersection form is played by the
Beauville-Bogomolov-Fujiki form, which we denote by
$q$ (see Subsection \ref{hk}). With this form, the
projectivization of the positive cone in $H^{1,1}_{\R}(X)$ is again viewed as a model of the
hyperbolic space. 

The classification of automorphisms of a hyperbolic
space plays here the same role as for complex surfaces:
elliptic automorphisms have finite order, 
loxodromic automorphisms have positive entropy
(\cite{_Yomdin:volume_}), and parabolic automorphisms, 
up to taking a power, preserve a
fibration in complex algebraic tori 
(this is Lo Bianco theorem,
\ref{_Parabolic_preserves_LoBianco_Theorem_}, \cite{LB}) 
and act ergodically in almost all of its fibers
(Theorem A). The existence of the fibration is in fact conditional on the following statement,
 known as the
{\bf hyperk\"ahler SYZ
  conjecture} (\cite{_Verbitsky:SYZ_}) or {\bf
  Lagrangian conjecture}, and verified for all known families of hyperk\"ahler manifolds (see \ref{ref-syz} for references).

\hfill

\conjecture\label{syz} Every nef Beauville-Bogomolov isotropic line bundle on an irreducible holomorphic symplectic manifold $M$ is semiample, that is, if $L$ is such a line bundle, then $L^{\otimes N}$ is base-point-free for some $N>0$.

\hfill

Under this conjecture,
the sections of  $L^{\otimes N}$, after eventually replacing $N$ by a multiple, define a fibration $\pi:
M\to X$ (Iitaka fibration associated to $L$, see \cite{Ueno}, chapter 5, and also \cite{Laz}, theorem 2.1.27, where the result is formulated for projective varieties but the proof does not use projectivity). By a famous result of Matsushita, $\pi$ is a
Lagrangian fibration, all smooth fibers of $\pi$ are tori
and the base is ``very similar'' to the projective space,
in particular it is a Fano variety of Picard number 1, possibly with mild singularities.
If it is smooth, it is biholomorphic to $\P^n$
(\cite{_Hwang:CP^n_}). 


Our main results can be stated as follows.

\hfill

{\bf Theorem A} (see \ref{_dense_or_main_Theorem_}): Let $M$ be an irreducible holomorphic symplectic manifold and $f:M\to M$ a parabolic automorphism preserving a Lagrangian fibration $\pi:M\to X$.
Then the orbits of $f$ on a sufficiently general
fiber of $\pi$ are dense in the euclidean 
topology. 

\hfill

{\bf Theorem B:} Let $G\subset \Aut(M)$ be a subgroup containing such an $f$. Then $G$ acts ergodically unless
$\pi$ is $G$-invariant, in other words, unless $G$ fixes the fixed point of $f$ on the absolute. This exception occurs if and only if $G$ is virtually abelian, or equivalently $G$ does not contain any loxodromic automorphisms.

\hfill

A slightly weaker formulation proved in the text is \ref{_ergodic_Theorem_}: see \ref{thB} which explains how to deduce Theorem B.

We refer to \cite{_Viana_Oliveira_}, \cite{M} for an introduction to ergodic theory and dynamics. Recall nevertheless the definition.

\hfill

A hyperk\"ahler manifold $M$ of maximal holonomy and of dimension $2n$ has a natural invariant probability measure $\mu$ coming from the volume form $\Omega^n\wedge\bar{\Omega^n}$, where $\Omega$ generates $H^{2,0}(M)$.

\hfill

\definition
Let $\Gamma$ be a group acting on a manifold $M$ and preserving a probability measure $\mu$.
We say that the action of $\Gamma$
is {\bf ergodic} if any $\Gamma$-invariant measurable subset
of $M$ is full measure or measure 0, with respect to $\mu$.

\hfill

\remark Equivalently, $(M, \Gamma)$ is ergodic
iff any $\Gamma$-invariant integrable function is constant almost everywhere.

\hfill

\remark\label{dense-ergodic-compact-group}
From ergodicity it follows that almost all
  orbits of $\Gamma$ are dense,
but the converse is not true
(\cite{_Furstenberg:torus_}). 
However if $g: K\to K$ is a translation of a compact group
$K$, the density of orbits implies the 
ergodicity (and in fact even a stronger property of unique
ergodicity, \cite{_Viana_Oliveira_}, section 6.3.3).



\hfill

Theorem A quickly implies Theorem B: when $\pi$ is not $G$-invariant, 
$G$ contains two parabolic automorphisms $p_1, p_2$
preserving distinct Lagrangian fibrations $\pi_1, \pi_2$, and we derive in the same way as in Cantat's paper
that the group generated by $p_1$ and $p_2$ acts
ergodically. 

\subsection{Existence of parabolic automorphisms}

Concerning the existence of parabolic automorphisms in the
presence of Beauville-Bogomolov isotropic nef line
bundles, the role of $(-2)$-curves 
on a K3 surface is played in higher dimension by MBM
classes (see \cite{AV-MBM} for their definition and
properties). If the Neron-Severi lattice of $M$ is of rank
at least three and contains no MBM
classes, then the automorphism group of $M$ is a lattice 
(``a lattice'' means in this context ``a discrete subgroup 
of finite covolume in $SO(\NS(M))$''),
as follows from the global Torelli
theorem. When $\NS(M)$ represents zero, any lattice in $SO(\NS(M))$ contains
a parabolic element; this
 implies the existence of a parabolic automorphism. 

To construct
maximal holonomy hyperk\"ahler 
manifolds with no MBM classes, we recall from \cite[Corollary 1.4]{AV-orbits} 
that the MBM classes have bounded square:
$0>q(z)>-N$, where $N$ depends only on the deformation class of $M$. Therefore to exhibit a lattice with a
parabolic automorphism we construct one
representing zero and not representing numbers between
zero and $-N$. A primitive embedding of such a lattice
into $H^2(M, \Z)$ gives a deformation of $M$ with exactly
this Neron-Severi lattice. 

Such a construction has been first carried out in
\cite{_AV:automorphisms_}, for loxodromic as well as
parabolic automorphisms. We give a particularly simple
parabolic version of these results in section
\ref{mbm-auto}: we observe that the absence of MBM classes
orthogonal to $L$, rather than the total absence, suffices
to produce parabolic automorphisms. This substantially
simplifies the construction (see also \cite{Mats-subgr} for another version of the argument).  

As indicated to us by the referee, a further study of dynamics of groups containing parabolic automorphisms of surfaces preserving distinct fibrations is carried out in \cite{CD}. It would be interesting to see whether it generalizes to hyperk\"ahler case.


\section{Basic notions}


\subsection{Hyperk\"ahler manifolds and Beauville-Bogomolov-Fujiki form}\label{hk}

For preliminaries on hyperk\"ahler manifolds, we refer to
\cite{_Besse:Einst_Manifo_} and
\cite{_Beauville_}. Throughout the paper, we denote by $M$
a complex manifold of the following kind.

\hfill

\definition
A {\bf maximal holonomy hyperk\"ahler manifold}, or an {\bf irreducible holomorphic symplectic manifold (IHS)} $M$ is a simply connected compact K\"ahler manifold 
such that $H^{2,0}(M)$ is generated by a nowhere degenerate form\footnote{By Calabi-Yau theorem, such an $M$ admits a unique hyperk\"ahler metric in each K\"ahler class, hence the names ``hyperk\"ahler'' and ``holomorphic symplectic'' are used interchangeably. See \cite{_Besse:Einst_Manifo_} for details.} $\Omega$.

\hfill

\theorem (Fujiki, \cite{_Fujiki_}).
Let $\eta\in H^2(M)$, and $\dim M=2n$, where $M$ is
hyperk\"ahler. Then $\int_M \eta^{2n}=cq(\eta,\eta)^n$,
for some primitive integer quadratic form $q$ on $H^2(M,\Z)$,
and $c>0$ a rational number.

\hfill

\definition This form is called
{\bf Bogomolov-Beauville-Fujiki form}. It is defined
by the Fujiki's relation uniquely if we also ask that the square of a K\"ahler class should be positive.

\hfill

\remark The form $q$ has signature $(3,b_2-3)$, and signature $(1, b_2-3)$ on the space of $(1,1)$-classes (see \cite{_Beauville_}) .
It is negative definite on primitive (1,1)-classes, and positive
definite on $\langle \Re \Omega, \Im\Omega, \omega\rangle$,
 where $\omega$ is any K\"ahler form. If $M$ is
 projective, $q$ has signature $(1, \rho-1)$
on the space $\NS(M)\otimes \R=H^{1,1}_\R(M)$, 
where $\NS(M)=H^{1,1}(M)\cap H^2(M,\Z)$ is 
the Neron-Severi group and $\rho=\rho(M)$  the Picard number.

\hfill

\remark\label{fuj-polar}
Recall that {\bf the polarization}
of a homogeneous polynomial map $P:\; W \arrow \Q$ 
of degree $d$ is a symmetric polynomial
$Q\in \Q[z_1, ..., z_d]$, $z_i \in W$
of degree 1 in each variable such that
$P(z)=Q(\underbrace{z, ..., z}\limits_{\text{$d$ times}})$.
Such a polynomial function $Q$ exists and is unique.
The polynomial $Q(z_1, ..., z_d)$ is proportional with
a positive rational coefficient
to the homogeneous component of 
$P\left(\sum_{i=1}^d z_i\right)$ which has degree 1
in each variable.
Applying this to the Fujiki relation
$\int_M \eta^{2n} = c q(\eta, \eta)^n$,
we obtain {\bf the polarized form} of Fujiki's relation,
used below in section \ref{AMGM-sect}:
$$\int_M \eta_1\wedge \dots \wedge \eta_{2n}=K\sum_{\sigma\in S_{2n}}q(\eta_{\sigma(1)},\eta_{\sigma(2)})\dots q(\eta_{\sigma(2n-1)},
\eta_{\sigma(2n)})$$ with $K=\frac{c}{(2n)!}$.

\subsection{MBM classes and existence of automorphisms}\label{mbm-auto}

Let $M$ be a projective maximal holonomy hyperk\"ahler
manifold. We denote by $\Pos_{\Q}(M)$\footnote{The notation $\Pos (M)$ usually refers to the positive cone in $H^{1,1}(M)$, so that $\Pos_{\Q}(M)$ is its ``algebraic version''.}  the connected component of the cone of classes
of positive square in $\NS(M)\otimes \R$ which contains the ample classes. The projectivization $\P(\Pos_{\Q}(M))$
shall be viewed as a model for the hyperbolic space $\H^k$, where $k=\rho-1$ and $\rho=rk(\NS(M))$ is the Picard number.

\hfill

\definition (\cite{Markman-survey}) 
The {\bf monodromy group $\Mon$} is the subgroup of
$GL(H^2(M,\Z))$ generated by the parallel transports
along the Gauss-Manin connection in families. The {\bf Hodge
monodromy group}, which we denote by $\Gamma$ here, is the
image in $GL(\NS(M))$ of the subgroup of 
$\Mon$ preserving the Hodge decomposition.

\hfill

\theorem (\cite{Verbitsky-Torelli}) The group $\Mon$ is a
finite index subgroup of the group 
$O(H^2(M,\Z), q)$, and $\Gamma$
is a finite index subgroup of $O(\NS(M), q)$. In
particular, $\Gamma$ acts on $\H^k = \P(\Pos_{\Q}(M))$
with finite covolume. \endproof

\hfill

\definition A class $z\in \NS(M)$ is MBM if for some
$\gamma\in \Gamma$, the orthogonal to $\gamma z$ contains
a face of the ample cone of some IHS birational model of
$M$. A class $z\in H^2(M, \Z)$ is MBM if it is MBM in some
complex structure where it is of type $(1,1)$.

\hfill

An MBM class $z$ always satisfies $q(z)<0$, and the dual
homology class to $z$, up to a scalar multiplier, is represented by a
rational curve on some deformation of $M$ (\cite{AV-MBM}).

\hfill


The following theorem is a restatement of results
proven in \cite{Markman-survey},
\cite{AV-hyperb} and \cite{AV-morkaw} (cf. \cite{AV-hyperb}, Theorem 3.8).

\hfill

\theorem\label{_Kahler_polyhedra_Theorem_}
Let $M$ be a projective hyperk\"ahler manifold,
and $\H_\Q:= \P(\Pos_{\Q}(M))$
the hyperbolic space associated with
its Picard lattice (which we assume to be of rank at least three). Denote by $S_\alpha\subset \H_\Q$
the hyperplane obtained as
orthogonal complement to an MBM class $\alpha$ of type (1,1), and consider the set $\{S_\alpha\}$ of all such hyperplanes.
Let $\Gamma\subset O(\NS(M))$ be the Hodge 
monodromy group of $M$. 
Then $\Gamma$ acts on the set
$\{S_\alpha\}$ with finitely many orbits. Moreover, denote by $K_1, ..., K_n$ the connected components of 
${\goth K}:=\frac{H_\Q\backslash \bigcup_\alpha S_\alpha}{\Gamma}$ with the orbifold structure given by the quotient map. 

Then
\begin{description}
\item[(i)]
For each $K_i$, the universal orbifold covering of $K_i$ is a 
locally finite convex hyperbolic polyhedron in
${\Bbb P}\Pos_{\Q}(M)$, obtained by projectivization of the
ample cone of a hyperk\"ahler birational model of $M$.

\item[(ii)] The set of isomorphism classes
of birational models of $M$ 
is in bijective correspondence with the set $\{K_i\}$.

\item[(iii)] Let $M_i$ be a birational model of $M$
such that $K_i$  is 
the image by the quotient map of the projectivization of its ample cone.
Denote by $\pi_1(K_i)$ the orbifold fundamental
group, that is, the deck group of the universal
cover of $K_i$. 
Then $\pi_1(K_i)=\Aut(M_i)/G$,
where $G$ is the group of automorphisms
of $M_i$ acting trivially on the Neron-Severi group.
Moreover, the group $G$ is finite.
\end{description}

{\bf Proof of \ref{_Kahler_polyhedra_Theorem_} (i):}
The set of orthogonal complements to 
MBM classes is locally finite, and
its $\Gamma$-quotient is finite, as shown in
\cite{AV-morkaw}. 
 In \cite{AV-MBM} it was proven that
the ample\footnote{More precisely, the results of \cite{AV-MBM} and \cite{Markman-survey} are formulated in the K\"ahler context, but
they carry over the ample setting without changes.} cones of birational models of $M$
are locally polyhedral, and that the orthogonal hyperplanes to the MBM classes
are precisely unions of codimension one faces of these cones.
It goes back to Markman \cite{Markman-survey}
that there is a wall-and-chamber decomposition of $\Pos_{\Q}(M)$ into the union of all $\Gamma$-images of ample 
cones of all birational models of $M$.
This takes care of \ref{_Kahler_polyhedra_Theorem_} (i).

\hfill

{\bf Proof of \ref{_Kahler_polyhedra_Theorem_} (ii):}
In \cite{Markman-survey} it was shown that the set
$U:=\frac{\text{The set of K\"ahler chambers in
    $\Pos(M)$}}{\Gamma}$ is in bijective correspondence 
with the set of isomorphism classes of 
birational models of $M$.  However, 
the set $U$ coincides by construction with the set
of connected components of ${\goth K}$.

\hfill

{\bf Proof of \ref{_Kahler_polyhedra_Theorem_} (iii):}
To see that $\pi_1(K_i)=\Aut(M_i)/G$,
let $\tilde K_i\subset H^{1,1}(M, \Q)$ be the ample cone of $M_i$.
Then $\tilde K_i \arrow K_i$ is the orbifold universal
covering. As shown in \cite{Markman-survey},
an element $u\in \Gamma$ is induced by
an automorphism of $M_i$ if and only if
$u$ maps its K\"ahler chamber to itself.
This gives $\pi_1(K_i)=\Aut(M_i)/G$.
To see that $G$ is finite
we use the Calabi-Yau theorem, which gives a 
bijection between the set of hyperk\"ahler
metrics on $M_i$ and the set of K\"ahler classes
in $H^{1,1}(M_i)$. If $u \in \Aut(M_i)$ fixes
a K\"ahler class $\omega \in H^{1,1}(M, \Q)$,
it preserves the corresponding hyperk\"ahler metric. Since $\Aut(M_i)$ is a Lie group by Bochner-Montgomery theorem, discrete in our case, and
the group of isometries of a compact manifold
is compact, this implies that $G$ is finite (compact and discrete).\footnote{This also follows from Fujiki-Liebermann theorem, see \ref{transl}}.
\endproof

\subsection{Cusps of a hyperbolic quotient of the K\"ahler
  cone and parabolic automorphisms}

\ref{_Kahler_polyhedra_Theorem_} associates an 
orbifold hyperbolic polyhedron $K$ with each K\"ahler chamber
of a hyperk\"ahler manifold, in such a way that
$\pi_1(K)$ is the automorphism group of the
manifold with this K\"ahler chamber.

In this subsection, we interpret the parabolic
automorphisms in terms of the cusp points of $K$,
which leads to a much simplified construction
of hyperk\"ahler manifolds admitting
parabolic automorphisms.

The following result is used to describe the cusps
geometrically. We use the ``thick-thin decomposition theorem'',
well-known in hyperbolic geometry and also valid in the orbifold context,
see \cite{Senska}.

\hfill

\theorem (Thick-thin decomposition)
 \\
Any $n$-dimensional complete hyperbolic 
manifold $\H/\Gamma$ of finite volume can be represented as
a union of a ``thick part'', which is a compact manifold
(with a boundary), and a ``thin part'', which is
a finite union of quotients of form $B/\Z^{n-1}$, where $B$ is a horoball
tangent to the boundary at a cusp point, and $\Z^{n-1} =\St_\Gamma(B)$.

\hfill

\proof See 
\cite[Section 5.10]{_Thurston:thick-thin_} or
\cite[page 491]{_Kapovich:Kleinian_}. \endproof

\hfill

A hyperbolic orbifold $\H/\Gamma$
is smooth if and only if $\Gamma$ is torsion-free.
By Selberg lemma 
(\cite{_Alperin:Selberg_}), any finitely-generated matrix subgroup
contains a finite index torsion-free subgroup. 
Applying this result to $\Gamma$, we obtain
that any hyperbolic orbifold $\H/\Gamma$ has a finite cover
which is a  hyperbolic manifold.

This gives the following version of thick-thin decomposition
theorem.

\hfill

\theorem (Thick-thin decomposition for orbifolds)
 \\
Consider $n$-dimensional complete hyperbolic 
orbifold $M:=\H/\Gamma$ of finite volume. Then $M$ can be represented as
a union of a ``thick part'', which is a compact manifold
(with a boundary), and a ``thin part'', which is
a finite union of quotients of form $B/\Gamma_i$, where $B$ is a horoball
tangent to the boundary at a cusp point, and 
$\Gamma_i =\St_\Gamma(B)$ contains $\Z^n$ as
a finite index subgroup. \endproof

\hfill

We are now in position to deduce the following corollary from
\ref{_Kahler_polyhedra_Theorem_}.

\hfill

\corollary\label{no_mbm} (see also 
\cite{Mats-subgr})\\ Let $M$ be a projective hyperk\"ahler
manifold with the Picard rank $\rho\geq 3$, carrying a nef isotropic line
bundle $L$. Assume that no MBM class is
Beauville-Bogomolov orthogonal to $L$. Then there exists an
automorphism of $M$ preserving $L$,  acting on the corresponding
hyperbolic space $\H= {\Bbb P}\Pos_{\Q}(M)$ parabolically.

\hfill

\proof  
Let $H:={\Bbb P}\Pos_{\Q}(M)/\Gamma$
be the hyperbolic orbifold associated 
with $M$, 
and $K\subset H$ a hyperbolic polyhedron obtained
as the image of the K\"ahler cone. As shown e. g. 
in \cite{_Kapovich:hyperbolic_}, the
cusps of $H$ are in bijective correspondence
with the set $R/\Gamma$ of orbits of $\Gamma$ on $R$, where
$R$ is the set of rational points  on the absolute.
The statement ``no MBM class is orthogonal to $L$''
is translated, in this language, as ``the hyperplanes
which cut $H$ into $K_1, ..., K_n$ don't pass through the
cusp associated with the isotropic vector 
$c_1(L)\in H^{1,1}(M,\Z)$''. This means
that the hyperbolic polyhedron $K$ 
contains an open neighbourhood of a cusp. 

Now, consider the thick-thin decomposition
of $K$ near the cusp. The map 
$\pi_1(B/\Gamma_i) \arrow \pi_1(K)$ 
is injective, because $\pi_1(K)$ acts
on $\tilde K\subset \H$ in a neighbourhood of a cusp,
and the quotient is isometric to $B/\Z^{n-1}$.
By construction, the group $\Gamma_i$
contains $\Z^{n-1}$ as a finite index subgroup.
All elements in $O(NS(M))\supset \Gamma$
corresponding to  $\Z^{n-1}\subset \Gamma_i\subset \Gamma$ 
are parabolic by construction. 
Since $\Aut(M)/G=\pi_1(K)$, this implies
that $\Aut(M)$ contains parabolic
automorphisms.
\endproof

\hfill

To finish this section, we give a simple method to
construct such hyperk\"ahler manifolds, inspired
by \cite{_AV:automorphisms_} but substantially simpler. The proof uses (also in \cite{_AV:automorphisms_}) the following result from \cite{AV-orbits}.

\hfill

\theorem (\cite{AV-orbits}, Corollary 1.4)
 Let $M$ be a hyperk\"ahler manifold with $b_2\geq 5$. There exists a positive integer $N$ depending only on the deformation class of $M$ such that 
 $q(z)\geq -N$ for any primitive MBM class $z$ in $H^2(M, \Z)$ (the smallest such number is called the MBM bound for $M$).

\hfill

\theorem Let $M$ be a hyperk\"ahler manifold with $b_2\geq 5$.
Then $M$ admits a deformation which has a parabolic automorphism.

\hfill

\proof When $b_2 \geq 5$, the
lattice $H^2(M,\Z)$ contains an isotropic vector by
Meyer's theorem \cite{Meyer}. Let $y$ be a primitive
element of $H^2(M,\Z)$ of square zero. Let $N$ be the MBM bound for $M$. 
We first show that $y$ can be embedded into a primitive
sublattice $\Lambda'$ of signature $(0,-)$\footnote{This notation means that $\Lambda'$ is of rank 2 as a $\Z$-module and that its quadratic form is negative semi-definite with 1-dimensional kernel}, not representing integers strictly between $-N$ and 0. 
Indeed one can find a rank-two
sublattice of signature $(+,-)$ of $H^2(M,\Z)$ orthogonal to $y$. Let $L$ be the smallest primitive sublattice containing that one and $y$: it has signature $(+,0,-)$. The quotient $\Lambda$ of $L$ by the sublattice spanned by $y$ is again of signature $(+,-)$.  
As such, $\Lambda$ has primitive vectors of arbitrarily negative
length\footnote{This elementary lemma can be checked as follows: if true for a lattice $\Lambda$, it is also true for $\Lambda '\subset \Lambda\otimes \Q$ commensurable 
to $\Lambda$; any $\Lambda$ is commensurable to a diagonal lattice, for which the statement is obvious.}. Take a primitive $u\in \Lambda$ with
square less than $-N$ and consider the sublattice $\Lambda'$ which is the inverse image of the one generated by $u$ under the projection map.
Then $\Lambda'$ is as required.

Now take any  $a\in H^2(M,\Z)$ with positive square and
not orthogonal to $y$.
By the surjectivity of the period map, 
there exists a deformation $M'$ of $M$ such that the smallest primitive sublattice containing $a$ and $\Lambda'$ is the Neron-Severi lattice of 
$M'$. 

Moreover, by \ref{_Kahler_polyhedra_Theorem_},
 the positive cone of $M'$ is
decomposed by the orthogonal hyperplanes to the MBM
classes into the union of K\"ahler chambers, that is, the K\"ahler cones
of the hyperk\"ahler birational models of $M'$ and their
monodromy images.
The hyperk\"ahler birational models of $M'$ are still
deformations of $M$ by \cite{Hu}; in fact such a
birational model is unseparable from $M'$ in the Teichm\"uller
space for hyperk\"ahler manifolds, which is not
Hausdorff. So replacing $M'$ by an appropriate birational
model, we may assume that 
$y$ is the class of a nef line bundle $L$ on $M'$.

Now the statement follows from \ref{no_mbm}. Indeed $NS(M')$ does not
contain MBM classes 
orthogonal to $y$, since all integer negative classes 
orthogonal to $y$ are in $\Lambda'$ (which is primitive by construction) and hence have square lower than $-N$.
\endproof

\subsection{Holomorphic Lagrangian fibrations on hyperk\"ahler manifolds}

\theorem (Matsushita, \cite{_Matsushita:fibred_})\\
Let $\pi:\; M \arrow X$ be a proper surjective holomorphic map with connected fibers
from a projective hyperk\"ahler manifold $M$ to a normal projective variety $X$, with $0<\dim X < \dim M$.
 Then $\dim X = \frac{1}{2} \dim M$, and the fibers of $\pi$ are 
Lagrangian (this means that the holomorphic symplectic
form vanishes on the fibers $\pi^{-1}(x)$). Moreover, $X$
has $\Q$-factorial log-terminal singularities and
Picard number 1. \endproof

\hfill

\definition A proper surjective holomorphic map $\pi:M\to X$ with connected fibers is called a fibration; if its fibers are Lagrangian, it is
{\bf a holomorphic Lagrangian fibration}. When $M$ is normal, $X$ may be assumed normal (taking the normalization). 
We assume it throughout the paper.

\hfill

\remark\label{fiber-torus-proj} It is well-known that a
smooth fiber $F$ of a Lagrangian fibration is a torus,
indeed it has trivial tangent bundle (the symplectic form
gives an isomorphism between the tangent and the conormal
bundle). Matsushita formulates his theorem in the projective setting but the proof is valid in the non-projective case as well.
When $M$ is hyperk\"ahler, $F$ is projective even
if $M$ itself is not (see Proposition 2.1 of \cite{Cam},
attributed to Voisin by the author; see also \cite[Corollary 3.4]{_SV:Moser_}). 
It is also true that $X$ is projective even if $M$ is not:
the reason is that $X$ is compact K\"ahler and Moishezon with rational
singularities, and such varieties are projective (see
\cite{KL}, Theorem 2.8, and 
\cite[Corollary 1.7]{_Namikava:projectivity_}).

\hfill

\remark The base $X$ of a Lagrangian fibration on a hyperk\"ahler
manifold is conjectured to be rational. Hwang (\cite{_Hwang:CP^n_}) proved that $X\cong \P^n$ whenever $X$ is smooth. Previously,
Matsushita (\cite{_Matsushita:higher_}) 
proved that $X$ has the same Dolbeault cohomology as
$\P^n$ if it is smooth.

\hfill







The inverse image of the ample generator of the Picard group of the base $X$ is clearly semiample with zero Beauville-Bogomolov square.
{\bf The hyperk\"ahler SYZ} \ref{syz}, when satisfied, implies the converse: any nef line bundle $L$ with $q(L)=0$
should be semiample, so that the linear system of sections of $L^{\otimes N}$ for suitable $N$ defines a Lagrangian fibration by Matsushita's theorem.

\hfill



\remark\label{ref-syz}
The hyperk\"ahler SYZ conjecture
has been proved for all known examples of hyperk\"ahler manifolds. See \cite{Bayer-Macri}
for Hilbert schemes of K3 surfaces, and implicitly \cite{Yoshioka} for
generalized Kummer varieties,
as well as \cite{Mongardi-Rap}, \cite{Mongardi-Onorati} for the two sporadic O'Grady's examples.

\hfill


\hfill

\proposition\label{transversal} Let 
$\pi_1:\; M \arrow X_1$,  $\pi_2:\; M \arrow X_2$
be two distinct Lagrangian fibrations.
Then the intersection of general fibers of 
$\pi_1, \pi_2$ is finite.
Moreover, the restriction of $\pi_2$ 
on any fiber of $\pi_1$ is surjective.

\hfill

\proof Let $\pi_1$, $\pi_2$ be two distinct Lagrangian fibrations given by the sections of nef line bundles $L_1$ resp. $L_2$. The Hodge index theorem 
implies $q(L_1, L_2)>0$. By Fujiki formula this implies that the intersection number
$L_1^nL_2^n>0$. This intersection number is also equal to the top self-intersection number 
of the restriction of $L_2$ to any fiber of $\pi_1$. The restriction
of $L_1$ to a smooth fiber of 
$\pi_2$ is therefore ample, because a nef and big line bundle on an
abelian variety is ample. Hence the restriction of $\pi_2$
to a smooth fiber of $\pi_1$ is a finite map.\footnote{Note however that $\pi_2$ can possibly contract an irreducible component of some singular fiber of $\pi_1$.} The
restriction to any fiber remains surjective, because it is
proper and dominant. \endproof


\section{Parabolic automorphisms of hyperkahler manifolds}


\subsection{Hyperk\"ahler manifolds with ergodic
  automorphism groups}




\hfill

\definition
An automorphism of a hyperk\"ahler manifold $M$
is called {\bf elliptic (parabolic, hyperbolic)}
if it is elliptic (parabolic, hyperbolic) 
as an element of $SO(H^{1,1}(M,\R))$ (cf. Theorem-definition 
from the Introduction and the remarks following it).

\hfill


\hfill

\remark
Let $p$ be a parabolic automorphism
of a hyperk\"ahler manifold, and $\eta$ be a fixed point
in $H^{1,1}(M)$ associated with the fixed point in the absolute.
Then $\eta$ is proportional to an integral cohomology class (as in the lines following the definition of a parabolic automorphism in the introduction)
which lies on the boundary of the K\"ahler cone, that is, the class of a nef line bundle $L$. Indeed the K\"ahler cone is invariant under the action of $p$, and $\eta$ can be obtained as a limit $p^i(w)$
for any K\"ahler class $w$ on $M$. 

\hfill

One of the main results of this paper can be stated as follows (see \ref{dense-ergodic} for the proof).

\hfill

\theorem\label{_ergodic_Theorem_}
Let $M$ be a hyperk\"ahler manifold
admitting two parabolic automorphisms
$p_1, p_2$ which preserve nef line bundles $L_1$, $L_2$,
with $c_1(L_1)$ not proportional to
$c_1(L_2)$.\footnote{Equivalently, the fixed points of $p_1$, $p_2$
on the absolute of the hyperbolic space ${\Bbb P}\Pos(M)$
are not equal.}
Assume the SYZ conjecture holds for $L_1$, $L_2$.
Then $p_1, p_2$ generate a group acting
on $M$ ergodically.

\hfill

\remark\label{thB}
If $M$ admits a parabolic automorphism $p$ and
$\Aut(M)$ is not virtually abelian,\footnote{A group is
  called {\bf virtually abelian} if it contains an abelian
  subgroup of finite index.} $M$ admits 
parabolic automorphisms
 which have different fixed points
on the absolute. Indeed it is well-known that a discrete subgroup in the stabilizer
of a point on the absolute 
in the group of isometries of a hyperbolic space is virtually abelian (see e.g. \cite{R}, chapter 5, Theorem 5.5.9). So if $\Aut(M)$ is not virtually abelian, then not all automorphisms of $M$ 
fix the same point $x$ of the absolute as $p$ does. Conjugating $p$ by an automorphism $g$ which does not fix $x$ we obtain
a parabolic automorphism $p'$ with a different fixed point $x'$. Note that $x'$ is also on the boundary of the K\"ahler cone,
since the latter is preserved by $g$. In particular, $x$ and $x'$ are classes of two different nef line bundles $L$, $L'$. Also, if 
SYZ conjecture holds for $L$, it holds for $L'$, since $L'$ is the image of $L$ by an automorphism.
Hence \ref{_ergodic_Theorem_} implies Theorem B from the introduction (if $\Aut(M)$ is virtually abelian and contains a parabolic automorphism,
all of its elements fix the same point $x$ on the absolute, as follows from \cite{R}, Theorem 5.5.9).

\hfill

\remark 
In \cite{_AV:automorphisms_} 
it was shown that that any hyperk\"ahler manifold
with $b_2 \geq 14$ admits a projective deformation $M_p$
with a parabolic automorphism; moreover, the automorphism
group of $M_p$ is arithmetic, hence $M_p$  admits 
parabolic automorphisms with different fixed points in the
absolute (note that the construction of \cite{_AV:automorphisms_} gives manifolds with Picard number $\rho\geq 5>2$). If we are looking for a single automorphism, rather than for two automorphisms with different fixed points on the absolute,
we now have a simpler construction in Subsection \ref{mbm-auto}.

\subsection{Parabolic automorphisms and Lagrangian
  fibrations on hyperk\"ahler manifolds}

Let $p: M\to M$ preserve a nef line
bundle $L$ defining a Lagrangian fibration
$\pi: M\to X$. Then $p$ takes fibers to fibers,
inducing an automorphism of the base $B$. It turns out
that this automorphism is of finite order.

This fact is due to Federico Lo
Bianco, see \cite[Theorem B]{_LoBianco:thesis_} for the case when $X=\P^n$, \cite{LB} for the general case.\footnote{Lo Bianco formulates his theorem in the projective context, but looking at the proof one sees that only the projectivity of the base $X$ is relevant, cf. \ref{fiber-torus-proj}.}

\hfill

\theorem \label{_Parabolic_preserves_LoBianco_Theorem_}
(Lo Bianco)
Let $p$ be a parabolic automorphism
of a hyperk\"ahler manifold $M$, and
$\pi:\; M \arrow X$ a Lagrangian
fibration such that for a K\"ahler class
$\omega$ on $X$, its pulback $\pi^*\omega$ is the
class on the boundary of the K\"ahler cone fixed by $p$.
 Then a certain positive iterate of $p$ preserves
the fibers of $\pi$.

\hfill

\remark
If $p$ satisfies $\pi\circ p=\pi$ (which by \ref{_Parabolic_preserves_LoBianco_Theorem_} is always the case after replacing $p$ by a power),  we say that $p$ {\bf preserves the Lagrangian 
fibration $\pi$.} The fibers (connected by definition of a 
Lagrangian fibration) are uniquely determined by $p$, because they are
uniquely determined by the cohomology class $\pi^*\omega$, and
$p$ fixes one and only one point on the absolute.

\hfill

\proposition\label{transl} Let $p$ be a parabolic automorphism preserving a Lagrangian fibration $\pi: M\to X$. Then some power of $p$ acts as a translation 
on a general fiber.

\hfill

\proof Let $T$ be a smooth fiber of $\pi$; 
this is a complex torus, algebraic by \ref{fiber-torus-proj} (or by \ref{_rank_1_Proposition_} below). 
Recall that an automorphism of a compact torus $A=\C^g/\Lambda$ is 
induced by an affine transformation of $\C^g$ with linear part preserving $\Lambda$, and $\Aut(A)^0$ is the group of translations, in general
not of finite index in $\Aut(A)$. If we denote by $\Aut^{\omega}(A)$ the subgroup of automorphisms preserving 
a K\"ahler class $\omega$, then by Fujiki-Liebermann theorem (see e.g. \cite{Fuj}) $\Aut(A)^0$ is of finite index in $\Aut^{\omega}(A)$. It turns out that the restriction of $p$ to $T$ must indeed
preserve a polarization. Indeed it is well-known
(\cite{_Oguiso:Picard_}, see also \ref{_rank_1_Proposition_}
here) that the rank of the restriction map $H^2(M,
\Z)\stackrel r \arrow H^2(T, \Z)$ is one. The map induced by $p$ must
preserve the image of $r$. This implies that the restriction 
of a K\"ahler class from $M$ to $T$
is also preserved by $p$, hence the conclusion. 
\endproof

\hfill


\remark
A complex torus $T$ is not a group, unless
you fix the origin. However, its translation group $\Aut(T)^0$ is a complex, commutative Lie group,
and it is isomorphic to $T$ as a manifold. We denote it by $T^0$ when we want to stress the group manifold structure.

\hfill

\claim
A translation $\tau_x$ of a torus $T$ by an element $x\in T^0$ 
has all its orbits dense if and only if $x$ is not contained in
a real subtorus $T'\subset T^0$. 
In this case, $\tau_x$ is ergodic 
(\ref{dense-ergodic-compact-group}).

\hfill

\subsection{Parabolic automorphisms are fiberwise ergodic on Lagrangian
  fibrations}

The following theorem is the main result of this paper,
used to obtain ergodicity of the parabolic actions.

\hfill

\theorem\label{_dense_or_main_Theorem_}
Let $p$ be a parabolic
automorphism of a hyperk\"ahler manifold
preserving a Lagrangian fibration $\pi:\; M \arrow X$.
Then there exists a full measure, Baire second category 
subset $R\subset X$, such that for all $r\in R$ the 
fibers $\pi^{-1}(r)$ are tori, and the automorphism 
$p$ acts on $\pi^{-1}(r)$ with dense orbits.

\hfill

\remark\label{dense-ergodic} \ref{_dense_or_main_Theorem_} immediately implies
\ref{_ergodic_Theorem_}, in the same way as for K3
surfaces in \cite{_Cantat:dyn_K3_TG_}.
Indeed, let $\Gamma$ be the group generated by
two parabolic automorphisms $p_1, p_2$, 
and $\pi_1, \pi_2:\; M \arrow X_i$ the Lagrangian fibrations
associated with $p_1, p_2$. Ergodicity of the action of
the group $\Gamma$ means that any $\Gamma$-invariant measurable
set $U$ is either full measure or zero measure. 
By Fubini theorem, the intersection of $U$
and the fibers $\pi^{-1}_i(x)$ is measurable with respect to Haar measure
for almost all $x\in X_i$ (where ``almost all'' is meant with respect to the pushforward of the canonical probability measure on $M$). Since $U$ is
$\Gamma$-invariant, and the parabolic
action on the fibers of $\pi_i$ is ergodic,
the intersection $\pi^{-1}_i(x)\cap U$ is full 
measure or zero measure for almost all $x\in X_i$.

Assume that $U$ is not of measure zero. 
Let $U_i\subset X_i$ be the set of all $x\in X_i$
such that $\pi_i^{-1}(x)\cap U$ is measurable and
of full measure in $\pi_i^{-1}(x)$. 
Applying Fubini theorem, we obtain that
the symmetric difference of $U$ and $\pi_1^{-1}(U_1)$, as well as that of 
$U$ and  $\pi_2^{-1}(U_2)$, is of measure zero. 

Hence the symmetric
difference of $\pi_2^{-1}(U_2)$ and $\pi_1^{-1}(U_1)$
is of measure zero, and applying Fubini again we see that this is only possible if both $\pi_2^{-1}(U_2)$ and $\pi_1^{-1}(U_1)$ are of full measure.
 Hence so is $U$.

\hfill

The proof of \ref{_dense_or_main_Theorem_}
is based on Hodge theory, and we give it in Section
\ref{_parabo_orbits_Section_} after the relevant results about variations
of Hodge structures are introduced.


\section{Variations of Hodge structures and Deligne's
  semisimplicity theorem}


We give a brief introduction to the variations of Hodge structures.
For more detail, see 
\cite{_Voisin-Hodge_,_Griffiths:transcendental_,_Peters_Steenbrink:MHS_}.

\subsection{Hodge structures}

\definition
Let $V_\R$ be a real vector space.
{\bf A (real) Hodge structure of weight $w$} 
on a vector space $V_\C=V_\R \otimes_\R \C$ 
is a decomposition $V_\C =\bigoplus_{p+q=w} V^{p,q}$, satisfying 
$\overline{V^{p,q}}= V^{q,p}$. It is called an {\bf integral} or {\bf rational}
Hodge structure if one fixes an integral or rational lattice $V_\Z$ or $V_\Q$
in $V_\R$. A Hodge structure is 
equipped with an $U(1)$-action\footnote{$U(1)$ denotes the unit circle in $\C$.} with $u\in U(1)$
acting as $u^{p-q}$ on $V^{p,q}$. {\bf  A morphism}
of integral/rational
Hodge structures is an integral/rational linear map which is $U(1)$-invariant.

\hfill

\definition 
{\bf A polarization}
on a rational Hodge structrure of weight $w$ is a $U(1)$-invariant 
non-degenerate 2-form $h\in V_\Q^*\otimes V^*_\Q$ 
(symmetric or antisymmetric depending on parity of $w$) which 
satisfies 
\begin{equation}\label{_Riemann_Hodge_Equation_}
	(-1)^{w(w-1)/2}\sqrt{-1}^{p-q}h(x, \bar x)>0 
\end{equation} 
(``Riemann-Hodge relations'')
for each non-zero $x\in V^{p,q}$.

\hfill

\definition
Two complex tori $T_1, T_2$ are called {\bf isogeneous} if
there exists a surjective finite holomorphic map $T_1 \arrow T_2$.

\hfill

\remark
The category of complex tori with a group structure (that is, the zero point fixed) is equivalent to the category of
integral Hodge structures of weight 1. 
The category of complex tori with a group structure up to isogeny 
is equivalent to the category of
rational Hodge structures of weight 1.
Under this correspondence, abelian varieties
correspond to Hodge structures admitting a polarization. 
See e.g. \cite[\S 7.2.2]{_Voisin-Hodge_}.

\hfill

\remark
The category $\cac$ of rational Hodge structures admitting 
a polarization is semisimple, that is, any object of $\cac$
is a direct sum of irreducible ones. In the case of weight one (the one relevant for this paper), 
an equivalent statement is that the category of 
abelian varieties up to isogeny is semisimple (Poincar\'e reducibility theorem).

\subsection{Variations of Hodge structures}

\definition
Let $X$ be a complex manifold.
A {\bf (polarized) real variation of Hodge structures (VHS)}
on $X$ is a complex vector bundle $(B, \nabla)$ with a flat connection
equipped with a parallel anti-complex involution (defining a fiberwise real structure) 
and (polarized) Hodge structures at each point, $B= \bigoplus_{p+q=w} B^{p,q}$
which satisfy the
following conditions: 
\begin{description}
\item[(a)] the polarization, the rational lattice 
and the real structure are preserved by $\nabla$.
\item[(b)]  (``Griffiths transversality condition'')
 $\nabla^{1,0}(B^{p,q}) \subset B^{p,q} \oplus
  B^{p+1,q-1}$.
\end{description}

\hfill

\example\label{family} (\cite{_Voisin-Hodge_}, 10.1, 10.2.1)
Let $\pi:\; M \arrow X$ be a proper holomorphic 
submersion. Consider the bundle
$V := R^k\pi_*(\C_M)$ with $k$-th cohomology
of $\pi^{-1}(x)$ as the fiber at $x$, the Hodge decomposition coming
from the complex structure on $\pi^{-1}(x)$,
and the Gauss-Manin connection. This defines
a variation of Hodge structures on $X$. Assume that there is a class $\omega$ in $H^2(M,\Z)$ such that its restriction on each fiber is ample. 
Then the cup product with $\omega$ induces a Lefschetz operator $L$, and the polarization is produced by setting $h(a,b)$ to be the product of $L^{n-k}a$ and $b$.In general one must restrict to the primitive cohomology to get the Hodge-Riemann condition. In this paper we deal with the special case $k=1$, so that the whole cohomology is primitive.

\hfill

\example\label{family-HK} More precisely, we shall consider the following situation. Take a lagrangian fibration $\pi: M\to X$ 
on a hyperk\"ahler manifold $M$. Let $X_0\subset X$ be the
locus of non-critical values of $\pi$ (by definition, the
critical locus includes the singularities of $X$,
hence $X_0$ is smooth). We denote by the same letter $\pi:
M_0\to X_0$ the restriction of the fibration to
$\pi^{-1}(X_0)$. This is a fibration in projective complex
tori. Then $R^1\pi_*(\C_{M_0})$ is a polarized variation
of Hodge structures even if $M$ is not projective. Indeed,
since by Oguiso's observation (also
\ref{_rank_1_Proposition_} here) the rank of the
restriction of $H^2(M,\Z)$ to the smooth fiber is 1, there
is a class $\omega$ in $H^2(M,\Z)$ with ample restriction
to the fibers.

\hfill

\claim
Polarized  integral variations of
Hodge structures of weight 1 are in a functorial bijection
with holomorphic Abelian fibrations. More precisely,
polarized variations of
Hodge structures of weight 1 correspond
to fibrations in complex tori $Y \arrow X$ with a holomorphic
section, giving the fiberwise choice of zero, 
and a line bundle $L$ on $Y$ which
is ample on all fibers, inducing a polarization. 
Given a fibration in abelian varieties without a section, one can
associate to it the ``relative Albanese'', or 
``translation fibration'', replacing each fiber $T$ by $T^0$. The variation of Hodge structures defined as in
 \ref{family} with $k=1$ is the same 
for a fibration and its relative Albanese.

\subsection{Deligne's semisimplicity theorem}

The following important theorem
is one of the main tools of our argument.

\hfill

\theorem \label{_Deligne_semisi_Theorem_}
(Deligne's semisimplicity theorem,
\cite{_Deligne:finitude_,_Wright:Deligne's_})\footnote{The theorem as stated by Deligne is valid for so-called complex variations of Hodge structures. This is actually a stronger statement which we shall not need in this paper. The version for real or rational VHS is obtained by putting together the summands conjugate by Galois action.} \\
Let $V$ be a polarized 
rational variation of Hodge structures over a quasiprojective base $M$.
Then the underlying flat bundle can be decomposed as
$V=\bigoplus_i W_i \otimes L_i$, where the
$L_i$ correspond to pairwise non-isomorphic
irreducible representations of $\pi_1(M)$, 
and $W_i$ are trivial representations. Moreover,
each $W_i$ is equipped with a Hodge structure, 
each $L_i$ is equipped with a variation of Hodge structures,
and the decomposition $V=\bigoplus_i W_i \otimes L_i$
is compatible with the Hodge structures on $W_i, L_i$.

\hfill

Applying this result in the case of weight 1, we notice
that one of two terms in $W_i \otimes L_i$
has weight 0, and the other has weight 1. Separating the sum in two parts according to whether $W_i$ or $L_i$ is of weight 0, we
obtain the following corollary. Recall that
a fibration in abelian varieties is called
{\bf isotrivial} if all its fibers are isomorphic.

\hfill

\corollary\label{product}
Let $V$ be a quasiprojective manifold fibered in abelian varieties. Then, after passing to an
isogeneous fibration $V_1$, we can decompose $V_1$
into a product of abelian fibrations with irreducible
monodromy of the Gauss-Manin connection and an isotrivial
abelian fibration. The isotrivial fibration corresponds to $\bigoplus_i W_i \otimes L_i$ with weight of $W_i$ equal to 1.
\endproof

\subsection{Deligne's global invariant cycle theorem}

Further on, we shall use the following important
theorem, also due to Deligne.

\hfill

\theorem \label{_inva_cycle_Theorem_}
(see for example \cite{_Voisin-Hodge_}, II, Theorem 4.24, also Proposition 4.23 and 4.18). 
Let $\pi: M_0 \arrow X_0$ be a smooth family of projective manifolds ,
and $V= R^i\pi_*(\C_{M_0})$ the corresponding Gauss-Manin
local system. Consider a monodromy invariant
vector $\alpha \in R^i\pi_*(\C_{M_0})\restrict x$,
and let $[\alpha]\in H^i(\pi^{-1}(x))$ be the corresponding
cohomology class in the fiber. Then $[\alpha]$
belongs to the image of the restriction map
$H^i(M) \arrow H^i(\pi^{-1}(x))$, where $M$ is any compactification of $M_0$.\footnote{$H^i(M)$ and $H^i(M_0)$ actually have the same image in $H^i(\pi^{-1}(x))$, \cite{_Voisin-Hodge_}, II, Proposition 4.23.}
\endproof

\hfill

\remark This theorem is usually formulated in the algebraic context, so that $M$ is projective and $M_0$ is quasiprojective. But the proof as given in \cite{_Voisin-Hodge_} adapts to our setting where $\pi: M\to X$ is a Lagrangian fibration on a hyperk\"ahler manifold, see 
\cite{_Voisin-Hodge_} II, remark 4.16 and compare to \ref{family-HK}.


\section{Variations of Hodge structures and Lagrangian
  fibrations}


In this section, we prove some preliminary results
about the Hodge structures associated
with a Lagrangian fibration.

\subsection{AM-GM inequality, products of Hermitian
  forms and Lagrangian fibrations}\label{AMGM-sect}

In this section we give an elementary proof of the following useful observation by Oguiso (\cite{_Oguiso:Picard_}, p.3)
who combined results by Voisin \cite{Voisin-lagrang} and Matsushita \cite{_Matsushita:higher_}.

\hfill

\proposition\label{_rank_1_Proposition_}
Let $\pi:\; M \arrow X$ be a Lagrangian fibration
on a hyperk\"ahler manifold and $T$ a smooth fiber\footnote{A projective torus, see \ref{fiber-torus-proj}} of $\pi$. 
Then the natural restriction map
$H^2(M)\arrow H^2(T)$ has rank 1.

\hfill

We use the following basic lemma, found in Cauchy's 
{\em ``Cours d'analyse de l'\'Ecole Royale Polytechnique,
premi\'ere partie, Analyse alg\'ebrique''} (1821).
In high school, it is usually called {\bf the
  AM-GM\footnote{Short for ``arithmetic mean - geometric
    mean''} inequality}.

\hfill 

\lemma Let $\alpha_1,..., \alpha_n$ be positive real numbers. Then
\[
\frac{\sum  \alpha_i}{n}\geq \sqrt[\stackrel{\text{\normalsize \em n}}{\;\;}]{\prod_i \alpha_i},
\]
and the equality holds if and only if all $\alpha_i$ are
equal.

\hfill

\corollary
Let $\alpha_1, ..., \alpha_n$ be positive numbers
such that $\sum\alpha_i=n$ and $\prod\alpha_i=1$.
Then all $\alpha_i=1$.
\endproof

\hfill

\proposition\label{_Hermitian_forms_basic_Proposition_}
Let $\omega_1, \omega_2$ be positive Hermitian forms
on a vector space $V=\C^{n}$, and $h_1, h_2$ the corresponding
Hermitian metrics.  Suppose that
$\omega_1 \wedge \omega_2^{n-1}= \omega_1^{n}=\omega_2^{n}$.
Then $\omega_1=\omega_2$.

\hfill

\proof Simultaneous diagonalization theorem implies 
that $h_1$ can be diagonalized in an orthonormal basis
for $h_2$. Let $A=\omega_1 \omega_2^{-1}$ be the corresponding
diagonal matrix.
Clearly, $\frac{\omega_1 \wedge \omega_2^{n-1}}{\omega_2^n}=\frac 1 n \Tr(A)$
and $\frac{\omega_1^{n}}{\omega_2^n}=\det A$. 
In other words, $\frac{\omega_1^{n}}{\omega_2^n}=\det A$
is the product $\prod_{i=1}^n \alpha_i$ of all eigenvalues of $A$, and
$\frac{\omega_1 \wedge \omega_2^{n-1}}{\omega_2^n}$
the arithmetic mean $\frac 1 n \sum_{i=1}^n \alpha_i$ of all these eigenvalues.
On the other hand, $\omega_1^{n}=\omega_2^{n}$
implies that the product of all eigenvalues of $A$
is 1, hence the geometric mean of these eigenvalues is
also 1. Then  
$\omega_1 \wedge \omega_2^{n-1}=\omega_1^{n}=\omega_2^{n}$  
implies that 
\[ \frac 1 n \sum_{i=1}^n \alpha_i=  \sqrt[n]{\prod_{i=1}^n
  \alpha_i}=1.
\] The AM-GM inequality implies that all $\alpha_i$ are
equal to 1, which gives $\omega_1=\omega_2$.
\endproof

\hfill


{\bf Proof of \ref{_rank_1_Proposition_}:}
Let $\eta_0\in H^2(X, \R)$ be a class in the singular cohomology
of $X$ satisfying $\eta_0^n=1$, where $n=\dim_\C X$,
$\omega_1, \omega_2\in H^2(M)$ any cohomology classes,
and $\eta:= \pi^* \eta_0$. We view $\eta$ as a de Rham cohomology class.
We need to show that $\omega_1 \restrict T$ is proportional
to $\omega_2 \restrict T$. Since $T$ is Lagrangian,
the (2,0)-forms are restricted to zero, and we may
assume that $\omega_1, \omega_2\in H^{1,1}(M)$.
Picking $\omega_1$ K\"ahler and replacing $\omega_2$
by an appropriate linear combination of $\omega_1$ 
and $\omega_2$, we may 
assume that $\omega_2$ is also K\"ahler.
Clearly, $q(\eta,\eta)=0$.
On the other hand, $\eta^n$ is cohomologous
to the fundamental class $[T]$ of the fiber of $\pi$,
because $ \eta_0^n=1$.
By Fujiki formula in the polarized form,
$\int_T \omega_i^n = \int_M \omega_i^n \wedge \eta^n=
Cq(\omega_i, \eta)^n$, where $C$ is the product of $K=c/(2n)!$ and of the number of permutations $\sigma\in S_{2n}$ such that for every pair of consequent indices 
$(2k-1, 2k)$ exactly one of the values $\sigma(2k-1)$ and $\sigma(2k)$ is greater than $n$.
Indeed apply \ref{fuj-polar} with $\eta_1=\dots=\eta_n=\omega_i$ and $\eta_{n+1}=\dots=\eta_{2n}=\eta$ and use $q(\eta,\eta)=0$.
In the same way and for the same $C$,
\[\int_T \omega_2\wedge\omega_1^{n-1} = 
\int_M \omega_2\wedge\omega_1^{n-1}\wedge \eta^n=
Cq(\omega_2, \eta) q(\omega_1, \eta)^{n-1}.
\]
Rescaling $\omega_i$ in such a way that
$q(\omega_1, \eta)=q(\omega_2, \eta)=C^{-1/n}$,
we obtain that 
\[ \omega_1 \wedge \omega_2^{n-1}\restrict T=
\omega_1^{n}\restrict T=\omega_2^{n}\restrict T
\]
for the two K\"ahler classes restricted to the fiber. 
Representing these K\"ahler classes by
translation-invariant Hermitian forms on $T$, we
obtain two Hermitian forms on $\C^n$
satisfying $\omega_1 \wedge \omega_2^{n-1}= \omega_1^{n}=\omega_2^{n}$.
By  \ref{_Hermitian_forms_basic_Proposition_},
these Hermitian forms are equal, hence  $\omega_1\restrict T$ 
is cohomologous to $\omega_2\restrict T$.
\endproof

\subsection{Irreducibility of the VHS associated with a Lagrangian fibration}

From now on, we denote by $X_0\subset X$ the locus of
non-critical values of $\pi:M\to X$, i.e. the maximal
smooth Zariski-open subset over which $\pi$ is smooth.

\hfill

\lemma\label{_parallel_2-forms_on_VHS_Lemma_}
Let $\pi:\; M \arrow X$ be a Lagrangian fibration
on a hyperk\"ahler manifold of maximal holonomy.
Let $V$ be the variation of Hodge structures over $X_0$ 
associated  with $R^1\pi_*(\C_M)$. Then 
the space of parallel sections of $\bigwedge^2 V$
is 1-dimensional.

\hfill

\proof
Let $\omega$ be the polarization on $V$ (see \ref{family-HK}).
 It is a parallel section of $\bigwedge^2 V^*$,
corresponding to the restriction of the Riemann-Hodge form
to the fibers of $\pi$. 
Since the smooth fibers of $\pi$ are tori,
the bundle of second cohomologies of the fibers of $\pi$ over $X_0$ is $\bigwedge^2 V$, as a variation of Hodge structures.
By Deligne's global invariant cycle theorem
(\ref{_inva_cycle_Theorem_}),
any monodromy invariant element in a fiber of  $R^2\pi_*(\C_M)$
is obtained as a restriction
of a globally defined cohomology class in $H^2(M)$.
Applying this to $\bigwedge^2 V=R^2\pi_*(\C_M)$,
we obtain that any antisymmetric, monodromy-invariant
2-tensor on $V$ is obtained as a restriction
of a globally defined cohomology class in $H^2(M)$.

By \ref{_rank_1_Proposition_},
the restriction map from $H^2(M)$ to a smooth fiber of $\pi$
has rank 1, hence the conclusion. 
Dually, any parallel 2-form on $V$
is equal to the polarization.
\endproof

\hfill


\hfill

\corollary \label{_mono_irre_Corollary_}
Let $\pi:\; M \arrow X$ be a Lagrangian fibration
on a hyperk\"ahler manifold. Then the Gauss-Manin local system $V:= R^1\pi_*(\C_M)$ over $X_0$
 is irreducible as a variation of Hodge
structures.  Moreover, the corresponding monodromy representation 
is irreducible, unless $\pi$ is isotrivial.

\hfill

\proof Suppose that $V$ is not irreducible as a rational 
variation of Hodge structures:
$V=V_1 \oplus V_2$. Restricting
the polarization $\omega$ to $V_1$ and $V_2$\footnote{The restriction of a polarization to a Hodge substructure is itself a polarization, in particular it is non-zero.},
we obtain that the space of parallel sections of
$\bigwedge^2 V^*\cong \bigwedge^2 V$ is at least 2-dimensional,
which is impossible by \ref{_parallel_2-forms_on_VHS_Lemma_}.

\hfill

To prove the second part, apply Deligne's semisimplicity  
\ref{_Deligne_semisi_Theorem_} to obtain
the decomposition $V=\bigoplus_i W_i \otimes L_i$.
Since $V$ is irreducible as a VHS, only one of the summands
is non-trivial: $V=W\otimes L$. Then either $W$ has
weight 0, that is, $W=\Q$, and
$L$ is irreducible of weight 1, 
or $W$ has weight 1 and $L$ has weight 0.
In the second case $\pi$ is isotrivial.
\endproof

\hfill

\remark More geometrically, one can argue that a Hodge substructure of $V$ yields a fibration in subtori on each smooth 
fiber of $\pi$. This in turn gives a rational map $\phi$ from $M$ to the component of the Chow variety of $M$ which parameterizes 
such subtori. Considering the inverse image by $\phi$ of an ample line bundle (or a K\"ahler class in the non-algebraic case, indeed cycle spaces 
on compact K\"ahler manifolds are compact and
K\"ahler by \cite{Varouchas}) on this component, one again sees that the rank of the restriction of $H^2(M)$ to a smooth fiber of $\pi$ must be greater than one: indeed this inverse image
is semiample but not ample on a smooth fiber of $\pi$.


\section{Parabolic automorphisms and their orbits}
\label{_parabo_orbits_Section_}


In this section and the next one, we prove \ref{_dense_or_main_Theorem_}
which claims that a generic orbit of a parabolic
automorphism of a hyperk\"ahler manifold
is dense in the corresponding fiber of a 
Lagrangian fibration.

\subsection{Parabolic automorphisms and monodromy-invariant
sub\-tori}


Let $p$ be an automorphism of 
a hyperk\"ahler manifold, preserving
a Lagrangian fibration $\pi:\; M \arrow X$ and acting as a translation on smooth fibers:
we recall that
for a regular value $x\in X$, 
the fiber $\pi^{-1}(x)$ is a compact complex (in fact algebraic) torus.
Consider the group of its translations
$\Aut(\pi^{-1}(x))^0= \frac{H^1(\pi^{-1}(x), \R)}{H^1(\pi^{-1}(x), \Z)}$ 
understood as the connected subgroup of its group of automorphisms.
Consider the element  $p_x\in \Aut(\pi^{-1}(x))^0$
induced by $p$, and let $P_x$ be the connected component of the unity of the closure of the group it
generates in  $\Aut(\pi^{-1}(x))^0$. This is a real subtorus.
Clearly the orbits of $p$ are dense in the fibers of $\pi$ if and only if $P_x$
is equal to $\Aut(\pi^{-1}(x))^0$, and otherwise their closures are finite unions of translates of subtori of the same dimension as that of $P_x$.

Since $\pi$ is differentiably locally trivial
over $X_0$, the 
bundle of all translations of the fibers of 
$\pi$ 
is equipped with a natural connection,
induced by the Gauss-Manin connection on $R^1\pi_*(\C_M)$.
Denote by $\Gamma$ its monodromy group.

\hfill

\proposition\label{_torus_constant_Proposition_}
In the setting described above, let
$r:=\max_{z\in X_0} \dim P_z$,
and assume that the dimension of  
$\dim P_x$ is equal to $r$. 
Then $P_x$ is a $\Gamma$-invariant
subtorus of $\Aut_0(\pi^{-1}(x))$.

\hfill

{\bf Proof:} The map
$x\mapsto p_x$ gives a section of
the relative Albanese fibration
of $\pi$. This fibration is locally
trivial in the real analytic category.
Then the map $x\mapsto p_x$ can be interpreted
as a map from a neighbourhood $U$ of some $x_0$ to a 
fixed real torus $T^n=\R^n /\Z^n$.

The dimension of $P_x$ can 
be expressed through the smallest number of
relations with integer coefficients between the coordinate components
of $p_x\in T_n$. For each of these integer relations,
denoted by $A$, the set $Y_A$ of $x$ for which 
the coordinate components of
$p_x$ satisfy the relation $A$ is real analytic in $U$.
Now let ${\mathcal P}$ be the set of all possible integer relations. For some $A\in {\mathcal P}$, 
 $Y_A=U$, so that the relation is satisfied everywhere on $U$. For other $A\in {\mathcal P}$ this set is a proper analytic subset.
Let $Y_U\subset U$ be the union of all $Y_A \subsetneq U$.
Since $Y_U$ is a countable union of measure zero sets,
its complement is of full measure, and 
outside of $Y_U$ the torus $P_x$ is constant (indeed it is defined by all relations $A$ such that $Y_A=U$).
Taking the union of $Y_U$'s we obtain a countable union of proper analytic subsets $Y\subset X_0$. The torus $P_x$ is locally constant of maximal 
dimension
outside $Y$, and over $Y$ some extra relations are
satisfied, so that $\dim P_y$ is lower for $y\in Y$.

Now by definition of the Gauss-Manin connection this means that $P_x$ is monodromy invariant. 
Indeed, a closed real analytic path  which is not contained in $Y$
intersects $Y$ in a countable set.
Clearly,  it suffices  to show that $P_x$ stays
constant along such a real analytic path.







Passing to the universal cover
of the fiber, we can consider $p$ as a real
analytic map from $\gamma$ to $\R^n$. 
Under this interpretation, $P_x$ is given
by the smallest rational subspace $V$ of $\R^n$
containing $p(x)$. This subspace is constant outside of $Y$, indeed it is defined by all integer relations $A$ satisfied over $X_0$, and more relations are satisfied over $Y$. This means that $p(\gamma) \subset V$, and $P_x= V$ 
for all $x\in \gamma\backslash Y$.
\endproof

\hfill

\remark\label{finite-order} A priori one could have $r=0$, however this means that $p$ is of finite order and therefore not parabolic. Indeed $r=0$ means that all translations in the fibers over $X_0$ are by torsion elements. But unless $p$ is of finite order, the fixed point set of $p^m$ is a proper analytic subvariety of $M$, and $M_0$ is not a countable union of proper analytic subvarieties.

\subsection{Ergodicity of parabolic automorphisms }

Now we can prove the following corollary,
which immediately implies 
\ref{_dense_or_main_Theorem_} when $\pi$ is not isotrivial.

\hfill

\corollary 
Let  $\pi:\; M \arrow X$  be a non-isotrivial
Lagrangian fibration on a hyperk\"ahler manifold
of maximal holonomy, and $p$ a parabolic automorphism
preserving $\pi$. Then $p$ acts with dense orbits on a general fiber of
$\pi$.

\hfill

\proof
The closure of an orbit of maximal dimension is given by 
a monodromy-invariant subspace in $H^1(\pi^{-1}(x))$, 
non-zero by \ref{finite-order}.
However, for a non-isotrivial Lagrangian fibration,
the monodromy is irreducible by \ref{_mono_irre_Corollary_}.
\endproof

\hfill



\section{Isotrivial Lagrangian fibrations}


Up to this point, to prove non-existence of subtori in a
holomorphic Lagrangian fibration,
we have used the irreducibility
of the corresponding Gauss-Manin local system.
This works when the Lagrangian fibration is 
non-isotrivial; for isotrivial Lagrangian
fibrations, an extra effort is needed. Indeed, though it
is  irreducible as a variation of Hodge structure, the
bundle of the first cohomologies is not necessarily
irreducible as a local system in this case (\ref{kummer}).

\hfill

The monodromy of the weight 1 variation of Hodge
structures $V$ associated with an isotrivial Lagrangian fibrations
is unitary, and therefore finite. 
Indeed, when the complex structure is constant in the family, the complex
structure operator $I$ acting on $V^{0,1}$ as $-\1$
and on  $V^{1,0}$ as $\1$ is parallel, and it
defines a positive definite Hermitian metric
\eqref{_Riemann_Hodge_Equation_}, which is also parallel, hence
preserved by the monodromy.

\hfill

\lemma\label{_torus_monodromy_Lemma_}
Let $V_\R$ be a real local system equipped with a lattice 
$V_\Z$, such that $V_\R = V_\Z \otimes_\Z \R$,
and suppose that the monodromy
of the torus bundle
$V_\R/V_\Z$ is trivial. Then the monodromy of
$V_\R$ is also trivial.

\hfill

\proof The group of 
automorphisms of a fiber of the lattice is isomorphic to
$GL(n,\Z)$, and this group acts faithfully on
the torus, indeed, already 
the action on its tangent space is faithful.
\endproof

\hfill

\lemma\label{_exist_section_from_para_Lemma_}
Let $\pi:\; M \arrow X$ be an isotrivial Lagrangian fibration
on a hyperk\"ahler manifold of maximal holonomy, $p$
a parabolic automorphism,
$M_0:= \pi^{-1}(X_0)$, and $P\subset M_0$ the family of subtori
obtained as components of the closure of orbits of $p$ 
(\ref{_torus_constant_Proposition_}).
Let $V:=R^1\pi_*(\C_{M_0})$ be the variation of Hodge structures
associated with $\pi$ and ${\Bbb P}\subset V_\R$
a parallel sub-bundle associated with $R^1\pi_*(\R_P)$.
Then ${\Bbb P}$ is locally generated by the
holomorphic sections of $V^{1,0}$ projected to $V_\R$
along $V^{0,1}$.\footnote{In other words, each fiber $P_x$ is generated, as a real vector space, by the images of local holomorphic sections at $x$.} 

\hfill

\proof $V^{1,0}/V_\Z \cong V_\R/V_\Z$ is the relative Albanese
space $\Alb_\pi$ of $\pi:\; M_0 \arrow X_0$.
Therefore, $p$ can be considered as a section
of the relative Albanese map
mapping $X_0$ to the total space of
${\Bbb P}/{\Bbb P}_\Z$. This section, considered
as a map $X_0 \arrow V^{1,0}/V_\Z$,
is holomorphic by definition. On the other hand,
the iterations of $p$ are dense in $P$,
hence after projection
generate ${\Bbb P}\subset V_\R$ .
\endproof

\hfill

The following proposition takes care of isotrivial
Lagrangian fibrations,
finishing the proof of \ref{_dense_or_main_Theorem_}.
Indeed, any parabolic automorphism whose orbits
are not dense in the fibers of the Lagrangian projection
defines a parallel subbundle ${\Bbb P} = (V_1)_\R \subset V_\R=R^1\pi_*(\R_{M_0})$
(\ref{_torus_constant_Proposition_}), and 
the corresponding subbundle in $V^{1,0}$ is generated
by holomorphic sections (\ref{_exist_section_from_para_Lemma_}).

\hfill

\proposition
Let $\pi:\; M \arrow X$ be an isotrivial Lagrangian fibration
on a hyperk\"ahler manifold of maximal holonomy.
Suppose that the corresponding
flat bundle $V:=R^1\pi_*(\C_{M_0})$ 
contains a parallel subbundle 
$V_1\neq 0$, a complexification of $(V_1)_\R\in V_\R$. If $(V_1)_\R$ is generated by real parts of
holomorphic sections of $V^{1,0}$, 
then $V_1=V$.

\hfill

\pstep
We prove that any proper parallel subbundle of $V$
is totally real (that is, contains no holomorphic subspaces) and Lagrangian with respect
to the polarization.

Since $\pi$ is isotrivial, the
complex structure operator on $V$ is parallel
with respect to the flat connection. 
Let $\omega$ be the polarization on $V$.
 It is a parallel section of $\bigwedge^2 V^*$,
induced by the restriction of the Riemann-Hodge form 
to the fibers of $\pi$.  Then the bilinear symmetric form 
$g(x, y)= \omega(I(x), y)$ is positive definite and parallel
on $V$.

Any parallel sub-bundle
$V_1 \subset V$  is Lagrangian with respect to
$\omega$. To see this, consider the orthogonal 
decomposition $V= V_1 \oplus V_1^\bot$ with respect to 
the positive definite form $g$ on $V$.
This decomposition is parallel, because $g$ is parallel.
Restricting $\omega$ to $V_1$, we obtain a parallel 2-form
on $V_1$; using the decomposition $V= V_1 \oplus V_1^\bot$,
we can consider this 2-form as a parallel section 
$\omega_1 \in \bigwedge^2(V)$. Since the
rank of $\omega_0$ is strictly less than 
$\rk \omega$, we have $\omega_1=0$ by 
\ref{_parallel_2-forms_on_VHS_Lemma_}. 
Therefore, $V_1$ is isotropic.
Similarly, the restriction of $\omega$
to $V_1^\bot$ also vanishes. This implies that
these two sub-bundles are Lagrangian. 

The fact that $V_1 \subset V$ is Lagrangian with respect to the polarization
implies that $V_1$ is totally real, indeed the polarization would be positive on a holomorphic subspace.

\hfill

{\bf Step 2:} 





Let $U \subset X$ be an open subset, and
$s:\; U \arrow V^{1,0}$ a holomorphic section
such that $\Re(s) \subset (V_1)_\R$.
Since $s$ is holomorphic, $\bar\6 s=0$
and $\nabla(s)=\nabla^{1,0}(s)$.
Unless $\nabla(s)=0$, this would imply that
$V_1\cap V^{1,0}\neq 0$; indeed,
$\nabla^{1,0}_\xi(s) \subset V^{1,0}$ for
any $\xi \in TX$. However, Step 1 implies
that $V_1\cap V^{1,0}=0$, implying that $\nabla (s)=0$.
Then $V_1$ is generated by parallel sections, which is
impossible by  Deligne's global invariant cycle theorem
(\ref{_inva_cycle_Theorem_}),
because $H^1(M)=0$ whenever $M$
is a hyperk\"ahler manifold of maximal holonomy.
\endproof

\hfill

\remark\label{kummer} Consider the Kummer K3 surface associated to a product of elliptic curves. The projection to each factor induces an isotrivial Lagrangian fibration on the Kummer surface, and
it is easy to see that the monodromy of its local system of the first 
cohomologies is $\{\pm 1\}$.
Hence this local system of real rank 2
is not irreducible. In the context of Deligne's theorem
\ref{_Deligne_semisi_Theorem_},
this example illustrates the case when
an irreducible real variation of Hodge structures
is obtained as a tensor product $L\otimes W$ of a rank 1
real local system $L$, understood as a variation
of Hodge structures of weight 0, and a weight 1 Hodge structure on a two-dimensional vector space $W$.

\hfill

\remark\label{singular} After this article has been released, Ljudmila Kamenova has remarked to us that similar results can also be obtained in the singular context, that is, for primitive symplectic varieties, where the analogues of the main theorems of hyperk\"ahler geometry have recently been obtained by Bakker and Lehn. See \cite{KL} for basic results on fibrations in this case (remarkably, the analogue of Matsushita's theorem holds, that is, any fibration is lagrangian with almost all of the fibers smooth projective tori, and projective base).

\hfill

{\bf Acknowledgements:} We are grateful 
to Serge Cantat, Ilyas Bairamov, Ljudmila Kamenova and Ben Bakker for interesting
discussions. Our special thanks go to the referee 
for his careful reading and pertinent remarks.

\hfill

{

{\small
\noindent {\sc Ekaterina Amerik\\
{\sc Laboratory of Algebraic Geometry,\\
National Research University HSE,\\
Department of Mathematics, 7 Vavilova Str. Moscow, Russia,}\\
\tt  Ekaterina.Amerik@gmail.com}, also: \\
{\sc Universit\'e Paris-11,\\
Laboratoire de Math\'ematiques,\\
Campus d'Orsay, B\^atiment 425, 91405 Orsay, France}

\hfill

\noindent {\sc Misha Verbitsky\\
{\sc Instituto Nacional de Matem\'atica Pura e
              Aplicada (IMPA) \\ Estrada Dona Castorina, 110\\
Jardim Bot\^anico, CEP 22460-320\\
Rio de Janeiro, RJ - Brasil\\
\tt  verbit@impa.br }\\
also:\\
{\sc Laboratory of Algebraic Geometry,\\
National Research University HSE,\\
Department of Mathematics, 7 Vavilova Str.\\ Moscow, Russia}.
 }
}}

\end{document}